\documentclass[11pt]{amsart}
\usepackage{amsmath,amsthm,amsfonts,amssymb,latexsym,epsfig}
\usepackage[T1]{fontenc}
\usepackage[latin1]{inputenc}
\usepackage{amsmath}
\usepackage{amssymb}
\usepackage{latexsym} 
\usepackage{amsthm}
\usepackage{amsfonts}
\usepackage{graphicx}
\usepackage{url}
\usepackage[right]{eurosym}
\usepackage{ifpdf}
\usepackage{array}

\numberwithin{equation}{section}
  
\theoremstyle{plain} 
\newtheorem{theorem}{Theorem}[section] 
\newtheorem{corollary}[theorem]{Corollary} 
\newtheorem{prop}{Proposition}[section]

\newtheorem{lemma}[theorem]{Lemma}

\theoremstyle{definition}

\title{$n$-LEVEL DENSITY OF THE LOW-LYING ZEROS OF QUADRATIC DIRICHLET $L$-FUNCTIONS}
\author{Peng Gao}
\address{Department of Computer and Mathematical Sciences,
University of Toronto at Scarborough, 1265 Military Trail, Toronto
Ontario, Canada M1C 1A4} \email{penggao@utsc.utoronto.ca}
\date{June 30 , 2008.}

\begin{document} 
\begin{abstract} 
   The Density Conjecture of Katz and Sarnak
associates a classical compact group to each reasonable
family of $L$-functions. Under the assumption of the Generalized Riemann Hypothesis, Rubinstein computed the $n$-level density of low-lying zeros for the family of quadratic Dirichlet $L$-functions in the case that the Fourier transform  $\hat{f}(u)$ of any test function $f$ is supported in the region $\sum^n_{j=1}u_j < 1$ and showed that the result agrees with the Density Conjecture.  In this paper, we improve Rubinstein's result on computing the $n$-level
density for the Fourier transform  $\hat{f}(u)$ being
supported in the region $\sum^n_{j=1}u_j < 2$.
\end{abstract} 

\maketitle

\tableofcontents 

\section{Introduction}
\label{sec1}

   The celebrated Riemann Hypothesis (RH)
   asserts that every non-trivial zero of the Riemann zeta function $\zeta(s)$ is of the form
   $\rho=1/2 + i\gamma$, with $\gamma$ real. 
   One possible approach towards establishing RH was suggested by P\'olya and
   Hilbert long time ago, who conjectured that one should find a
   Hilbert space and a suitable linear operator acting on it whose
   spectrum is given by the non-trivial zeros of $\zeta(s)$. Then
   depending on the properties of the operator, the
   zeros are forced to lie on a line.

   Montgomery \cite{Mont} initiated the study of the local spacings
   of the zeros of $\zeta(s)$. Assuming RH, he derived, under certain
   restrictions, the pair correlation of the zeros of $\zeta(s)$ and
   made the famous pair correlation conjecture. His result drew the attention of
   Dyson, who pointed out that the eigenvalues of Gaussian Unitary Ensemble
(GUE) have the same pair correlation as $N \rightarrow
\infty$. With this insight, Montgomery went on to conjecture (the
GUE Conjecture) that all the correlations would match up for the zeros of $\zeta(s)$
and eigenvalues of Hermitian matrices. Montgomery's pair correlation result then
provides the first evidence in favor of the spectral interpretation
of the zeros, suggesting that the relevant operator might be
Hermitian.

  Rudnick and Sarnak
\cite{R&S} studied the distribution of zeros of
 general $L$-functions.
 According to the Langlands program (see \cite{Ge} , \cite{Murty}, \cite{knapp}),
  all $L$-functions can be written as products of ``standard''
 $L$-functions $L(s, \pi)$
attached to automorphic cuspidal representations of $GL_M$ over
${\mathbb Q}$. The Generalized Riemann
Hypothesis (GRH) for $L(s, \pi)$ asserts that the non-trivial zeros
of $L(s, \pi)$ all lie on the line $\Re s =1/2$ when properly normalized. The result of Rudnick and
Sarnak \cite{R&S} asserts that the $n \geq 2$ correlations of the
zeros of $L(s, \pi)$ (assuming GRH and with some restrictions) are
universally the GUE ones. Their result is confirmed by the numerical
work of Rumely \cite{Rum} for Dirichlet $L$-functions and of
Rubinstein \cite{Ru} for a variety of $GL_2/{\mathbb Q}$
$L$-functions.

 On the random matrix theory side, Katz and Sarnak showed
  \cite{K&S1} that the local spacing statistics for the following
  types of spaces are universally GUE: (special) orthogonal $SO(\text{odd}), SO(\text{even})$,  (special) unitary $SU(N)$ and
   (unitary) symplectic $USp(2N)$ groups.  In contrast to this, they also showed
  that the distributions of the eigenvalues nearest to
  $1$ is sensitive to the particular symmetry of the group. Based on this and studies on the function field
    analogues, Katz and Sarnak conjecture (the Density Conjecture, \cite[p. 20]{K&S}) that to
  each reasonable family $\mathcal{F}$ of $L$-functions (by which we mean a collection of geometric
objects and their associated $L$-functions, where the geometric
objects have similar properties), one
  can associate a classical compact group $\mathcal{G}(\mathcal {F})$ mentioned above,
  such that the behavior
of zeros near the central point $s = 1/2$ (the low-lying zeros) of
$L$-functions in $\mathcal{F}$ is the same as the behavior of
eigenvalues near $1$ of matrices in $\mathcal{G}(\mathcal {F})$. This
conjecture further suggests that in terms of the spectral
interpretation, the relevant operators for $L$-functions should have
some symmetry corresponding to the relevant classical group.

  The  $1$-level density for quadratic twists of
$\zeta(s)$ was obtained by \" Ozl\" uk and Snyder \cite{O&S}. A stronger result which applies for $\zeta(s)$ as well as all $L(s, \pi)$ was obtained
by N. Katz and Sarnak \cite{K&S}. The
general case $n \geq 1$ was obtained by Rubinstein \cite{R}. Results are also obtained for a wide variety of
families of $L$-functions, see \cite{K&S1}, \cite{ILS},
\cite{Miller}, \cite{Gu}, \cite{D&M} and \cite{H&M}. In all the results above, the Fourier transforms for
test functions are assumed to be compactly supported in certain
regions. In this paper, we will study the $n$-level density of the
low-lying zeros of quadratic Dirichlet $L$-functions. This family of $L$-functions was studied by Rubinstein in \cite{R}. Our result in this paper
extends the one of Rubinstein's by enlarging the region in which the
Fourier transform for the test function is supported. In the next section, we shall give a detailed description of the problem studied in this paper and give an explanation on the method being used in the proof of the main result. Section \ref{sec3.0}  will be devoted entirely for the proof of our main result of the paper.

\section{Statement of the Result and Outline of the proof}
\label{sec2.1}
    In this paper, we consider the average
behavior of the low-lying non-trivial zeros of a family of
$L$-functions hoping to find evidence in favor of the Density
Conjecture. We focus on the quadratic twists of $\zeta(s)$, $\{L(s,
\chi_{d} ) \}$ as our family of $L$-functions, where $\chi_d(n) =
(\frac {d}{n})$ is the Jacobi symbol. We recall
    here  all real (quadratic) non-principal characters are of the form
    $\chi_d(n)$, where $d$ belongs to the
    set $Q$ of quadratic discriminants, $Q = \{ d: d $ is not a square and
    $d \equiv 0 \hspace{0.1in} \text{or} \hspace{0.05in} 1 \hspace{0.05in} (\text{mod} \hspace{0.05in} 4) \}$.
    The character $\chi_{d}(n)$
    is induced by a primitive character $\chi_{d'}(n)$ where $d'$
    belongs to the set $F$ of fundamental discriminants,
    $F = \{ d: d \equiv 1 \hspace{0.05in} (\text{mod} \hspace{0.05in} 4), d$ square-free $\}$ $\cup$
    $ \{ d: 4|d, d/4 \equiv 2 \hspace{0.05in} \text{or} \hspace{0.05in} 3
    \hspace{0.05in} (\text{mod} \hspace{0.05in} 4), d/4$ square-free
    $\}$ (see \cite[\S 5]{D}).

  To facilitate our study, it is more convenient to consider the family of $\{L(s, \chi_{8d}
) \}$ with odd, positive square-free $d$.  We may thus avoid the
possible impact of $\chi_2$ in our study and it follows from the
discussions in the proceeding paragraph that $\chi_{8d}$ is a real,
primitive character with conductor $8d$ and $\chi_{8d}(-1)=1$. Let
$D$ denote the set of such $d$'s, and let
    $D(X)=\{d \in D: X \leq d \leq 2X \}$. It's easy to see that the
    cardinality $|D(X)|$ of $D(X)$ is asymptotically $4X/\pi^2$.
    In
this case, the Density Conjecture suggests a unitary symplectic
symmetry  and our main result below (Theorem
\ref{thm2.1} ) or its corollary provides an evidence for it.

  Let $h$ be a function on ${\mathbb R}$, we say
that $h$ tends to
    $0$ rapidly at infinity (or rapidly decreasing) if for each positive
    integer $m$, the function $x \mapsto (1+|x|)^mh(x)$ is
    bounded for $|x|$ sufficiently large. We define the Schwartz
    space $S({\mathbb R})$ to be the set of functions on ${\mathbb
    R}$ which are infinitely differentiable, and which tend to $0$
    rapidly at infinity, as well as their derivatives of all orders (see \cite[Chap. VIII, \S 4]{L}).

    Now for each $f_i$ even and in $S({\mathbb R})$, we write
\begin{eqnarray*}
  f(x_1, \ldots, x_n) &=& \prod^{n}_{i=1}f_i(x_i),  \\
  \hat{f}(u_1, \ldots, u_n) &=& \prod^n_{i=1}\hat{f}_i(u_i).
\end{eqnarray*}

  Throughout the paper, we assume GRH and
    write the non-trivial zeros of $L(s, \chi_{8d})$ as
\begin{equation*}
\label{01.0}
    \frac 1{2}+i\gamma^{(j)}_{8d}, j=\pm 1, \pm 2, \ldots,
\end{equation*}
    where
\begin{equation*}
   0 \leq \gamma^{(1)}_{8d} \leq \gamma^{(2)}_{8d} \leq
   \gamma^{(3)}_{8d} \ldots
\end{equation*}
    and 
\begin{equation}
\label{01.00}
    \gamma^{(-j)}_{8d}=-\gamma^{(j)}_{8d}.
\end{equation}

     Let $N(T, 8d)$ denote the number of $\gamma^{(j)}_{8d}$'s with
     $|\gamma^{(j)}_{8d}|<T$,  then we have (see \cite[\S 16]{D}) for $T \geq
     2$
\begin{equation*}
   N(T, 8d)=\frac {T}{\pi}\log (\frac{8dT}{2\pi})-\frac
   {T}{\pi}+O \left (\log T+ \log \left(8d \right ) \right).
\end{equation*}
   We note here the above formula doesn't provide us much
     information when $T$ is close to $1$, since then the main term
     and the error term are of the same size. However, by taking $T=1$ in
     $N(T, 8d)$ one may think that there are about $\log (8d)/\pi$ many zeros up to height $1$ so that the mean
     spacing between them is $2\pi/\log(8d)$, which is also likely to be the
     height of the lowest zero. Thus by scaling back this height to
     $1$, and note that $d$ is of size $X$, we are then
     interested in studying the asymptotic behavior of
\begin{equation*}
     \lim_{X \rightarrow \infty} \frac {\pi^2 }{4 X}\sum_{d \in D(X)}{\sum_{j_1, \ldots,
     j_n}}^*f(L\gamma^{(j_1)}_{8d}, L\gamma^{(j_2)}_{8d}, \ldots, L\gamma^{(j_n)}_{8d})
\end{equation*}
     with $L=\log X/(2\pi)$ and
   $\sum^*_{j_1, \cdots,
     j_n}$ is over $j_k=\pm 1, \pm 2, \ldots, $ with $j_{k_1} \neq \pm j_{k_2}$ if $k_1 \neq
     k_2$. Note that since $f$ is rapidly decreasing, only the low-lying zeros
     contribute to the sums above.
     The Density Conjecture predicts that
\begin{equation}
\label{01.01}
    \lim_{X \rightarrow \infty} \frac {\pi^2 }{4 X}\sum_{d \in D(X)}{\sum_{j_1, \ldots,
     j_n}}^*f(L\gamma^{(j_1)}_{8d}, L\gamma^{(j_2)}_{8d}, \ldots, L\gamma^{(j_n)}_{8d})
    =\int_{{\mathbb R}^n}f(x)W^{(n)}_{\text{USp}}(x)dx,
\end{equation}
     where
\begin{equation*}
   W^{(n)}_{\text{USp}}(x_1, \ldots, x_n)=\det(K_{-1}(x_j,x_k))_{\substack{1\leq j \leq n \\ 1 \leq k \leq
   n}},  
\end{equation*}
   and
\begin{equation*}
    K_{-1}(x, y)   = \frac {\sin (\pi (x-y))}{x-y}-\frac {\sin (\pi (x+y))}{x+y}.
\end{equation*}   

   Rubinstein \cite[Theorem 3.1]{R} has shown that
\eqref{01.01} holds for $\hat{f}(u_1, \ldots, u_n)$ being supported in the region
$\sum^n_{i=1}|u_i|<1$. Presumably, \eqref{01.01} also holds for $f$
being compactly supported (by the Density Conjecture). However,
since $f$ and $\hat{f}$ can't both be compactly supported, we can't
approximate a compactly supported function by functions with
compactly supported Fourier transforms. Thus one needs to prove
\eqref{01.01} for the support of $\hat{f}$ being as large as
possible. The main result in this paper is an evaluation of the
left-hand side of \eqref{01.01} for $\hat{f}$ supported in
$\sum^n_{i=1}|u_i|<2$, which is given in the following:
\begin{theorem}
\label{thm2.1}
   Assume GRH and assume in \eqref{01.01} that each $f_i$ is even and in $S({\mathbb R})$ and $\hat{f}$ is supported
   in $\sum^n_{i=1}|u_i|<2$. Then
\begin{eqnarray}
\label{10.06}
   && \lim_{X \rightarrow \infty} \frac {\pi^2 }{4 X}\sum_{d \in D(X)}{\sum_{j_1, \ldots,
     j_n}}^*f(L\gamma^{(j_1)}_{8d}, L\gamma^{(j_2)}_{8d}, \ldots,
     L\gamma^{(j_n)}_{8d}) \\
  &=& \sum_{\underline{F}}(-2)^{n-\nu(\underline{F})}
\left ( \prod^{\nu(\underline{F})}_{l=1} \left ( |F_l|-1 \right)!
\right )\sum_{S}\left ( \prod_{l \in S^c}
\int_{\mathbb{R}}F_l \left(x \right )dx \right )   \nonumber \\
 && \cdot  \sum_{S_2 \subseteq S}\left ( \left (\frac {-1}{2} \right )^{|S^c_2|}
 \prod_{l \in S^c_2} \int_{{\mathbb R}}\hat{F}_l(u)du \right )
  \nonumber \\
&&  \cdot \left ( \left(\frac {1+(-1)^{|S_2|}}{2}\right )2^{\frac
{|S_2|}{2}}\sum_{(A;B)}\prod^{|S_2|/2}_{i=1}
\int_{\mathbb{R}} |u_i| \hat{F}_{a_i}(u_i)\hat{F}_{b_i}(u_i)du_i  \right. \nonumber \\
 &&- \frac {1}{2}\sum_{\substack {S_3 \subsetneq S_2 \\ |S_3| \text { even}}}2^{\frac
{|S_3|}{2}}
 \left ( \sum_{(C;D)}\prod^{|S_3|/2}_{i=1}\int_{\mathbb{R}} |u_i| \hat{F}_{c_i}\left (u_i \right )
 \hat{F}_{d_i}\left(u_i \right )du_i \right )
  \nonumber \\
 && \cdot  \left.  \sum_{I \subsetneq S^c_3}(-1)^{|I|}(-2)^{|S^c_3|}
 \prod_{i \in I}\int^{\infty}_{0} \hat{F}_i(u_i)du_i
 \int_{\substack{(\mathbb{R}_{\geq 0})^{I^c} \\ \sum_{i \in I}u_i \leq \sum_{i \in I^c}u_i-1}}\prod_{i \in
I^c}\hat{F}_i(u_i)du_i  \right ), \nonumber
\end{eqnarray}
  where $\underline{F}$ ranges over all ways of decomposing $\{1, \ldots, n \}$ into disjoint subsets
$\{ F_1, \ldots,$ $ F_{\nu(\underline{F})} \}$ with
\begin{equation}
\label{8.02}
   F_l(x)=\prod_{i \in F_l}f_i(x),
\end{equation}
  $S$ ranges over all $2^{\nu(\underline{F})}$ subsets of $\{1, 2,
\ldots, \nu(\underline{F}) \}$ and $S^c$ denotes the complement of
$S$. Here $\sum_{(A;B)}$ is over all ways of pairing up the elements
of $S_2$ and $\sum_{(C;D)}$ is over all ways of pairing up the
elements of $S_3$. Empty products are taken to be $1$.
\end{theorem}

  Unfortunately, we are not able to show the above result matches with the right-hand side
  of \eqref{01.01} in general. In the case of $n \leq 3$, one can show by direct computation that our result above does match with the right-hand side of \eqref{01.01}. As the proof is not particularly enlightening, we shall omit it here by stating the following result and refer the interested reader to \cite{Gao} for the proof:
\begin{corollary}
\label{cor2.01}
   Assume GRH. For $n \leq 3$, equality holds in \eqref{01.01} if each $f_i$ is even and in $S({\mathbb R})$ and $\hat{f_i}$ is supported
   in $\sum^n_{i=1}|u_i|< 2$.
\end{corollary}
  We note here one can also check that Theorem \ref{thm2.1} implies Rubinstein's result.

\subsection{Outline of the Proof of Theorem \ref{thm2.1}}
\label{sec 2.3}
     Our approach of proving Theorem \ref{thm2.1} relies on the
     following explicit formula, which essentially converts the sum
     over zeros of the $L$-function to the sum over primes, the one we
     use here follows from \cite[(2.16)]{R&S}:     
\begin{lemma}
\label{lem2.2.1}
   Let $h$ be a function on $\mathbb{R}$ with $\hat{h}$ smooth and compactly supported. Then
\begin{equation}
\label{01.3}
   \sum_{j}h(L\gamma^{(j)}_{8d})=\int^{\infty}_{-\infty}h(t)dt-\frac 1{2}
   \int^{\infty}_{-\infty}\hat{h}(u)du-S(d,X; \hat{h})+O(\frac {\log \log X}{\log
   X}),
\end{equation}
    with the implicit constant depending on $h$ and
\begin{equation*}
   S(d,X;\hat{h})=\frac 1{\log X}\sum_p\frac {\log
   p}{\sqrt{p}}\genfrac(){1pt}{}{8d}{p}\left ( \hat{h}(\frac {\log p}{\log X})+\hat{h}(\frac {-\log p}{\log
   X})\right
   ).
\end{equation*}
\end{lemma}

  Note that $f_i(-x)=f_i(x)$ implies that
\begin{eqnarray}
\label{8.1}
 && \frac {\pi^2 }{4 X}\sum_{d \in D(X)}{\sum_{j_1, \ldots,
     j_n}}^*f(L\gamma^{(j_1)}_{8d}, L\gamma^{(j_2)}_{8d}, \ldots, L\gamma^{(j_n)}_{8d})  \\
 &=& \frac {2^n \pi^2 }{4 X}\sum_{d \in D(X)}\sum_{\substack{j_1, \cdots,
     j_n \\ \text{positive} \\ \text{and} \\ \text{distinct} }}\tilde{f}_{8d}(j_1, \ldots,
     j_n), \nonumber
\end{eqnarray}
   where
\begin{equation*}
  \tilde{f}_{8d}(j_1, \ldots, j_n)=\prod^n_{i=1}f_i(L\gamma^{(j_i)}_{8d}).
\end{equation*}
   The sum in \eqref{8.1}  is over distinct indices $j_1, \ldots,
   j_n$. In order to apply Lemma \ref{lem2.2.1} to \eqref{8.1}, we
   would prefer to have a sum over all indices. To remove the distinctive
   condition, we use a combinatorial sieving as in \cite[p.
   305]{R&S}. We begin with some set-theoretic combinatorics. A set
   partition $\underline{F}$ of $\underline{N}=\{1, \ldots, n \}$ is
   a decomposition of $\underline{N}$ into disjoint subsets $\{ F_1, \ldots, F_{\nu(\underline{F})}
   \}$. The collection $\Pi_n$ of all set partitions of
   $\underline{N}$ forms a lattice with the partial ordering given
   by $\underline{F} \preceq \underline{G}$ if every subset $G_i$ of $\underline{G}$ is a union
   of subsets of $\underline{F}$. The minimal element of $\Pi_n$ is
   $\underline{O}=\{ \{1\}, \{2\}, \ldots, \{n\}\}$ and the maximal
   element is $\underline{N}=\{ 1, 2, \ldots, n \}$.

   The M\"obius function of a poset such as $\Pi_n$ is the unique
   function $\mu(x,y)$ so that for any functions $g,h: \Pi_n
   \rightarrow {\mathbb R}$, satisfying
\begin{equation*}
   g(x)=\sum_{x \preceq y}h(y),
\end{equation*}
   we have
\begin{equation*}
   h(x)=\sum_{x \preceq y}\mu(x,y)g(y).
\end{equation*}
   In the case of $\Pi_n$, the M\"obius function can be computed
   explicitly \cite[\S 25]{L&W}, in particular
\begin{equation*}
   \mu(\underline{O},
   \underline{F})=\prod^{\nu(\underline{F})}_{j=1}(-1)^{|F_j|-1}(|F_j|-1)!.
\end{equation*}

   Given a set partition $\underline{F}=\{ F_1, \ldots, F_{\nu(\underline{F})}\}\in
   \Pi_n$, we have an embedding $l_{\underline{F}}: {\mathbb R}^{\nu} \rightarrow {\mathbb R}^n$,
$l_{\underline{F}}(x_1, \ldots, x_{\nu(\underline{F})})=(y_1,
\ldots, y_n)$ with $y_i=x_j$ if $i \in \underline{F}_j$. For
example, for $n=3$, the possible $\underline{F}$'s  are $\{ \{1, 2,
3 \} \}$, $\{ \{1, 2 \}, \{ 3 \} \}$, $\{ \{1, 3 \}, \{ 2 \} \}$,
$\{ \{ 1 \}, \{2, 3 \} \}$, $\{ \{1 \}, \{ 2 \} , \{ 3 \} \}$ and
$l_{\{ \{ 1 \}, \{2,3 \} \}}(x_1, x_2)=(x_1, x_2, x_2)$.

   We now define
\begin{eqnarray*}
   R_{{\underline{F}}}(f) &=& \sum_{\substack{j_1, \ldots,
     j_n \\ \text{positive} \\ \text{and} \\ \text{distinct} }}
     \tilde{f}_{8d}(l_{\underline{F}}(j_1, \ldots,
     j_{\nu(\underline{F})})), \\
   C_{{\underline{F}}}(f) &=& \sum_{\substack{j_1, \ldots,
     j_n \\ \text{positive}}}
     \tilde{f}_{8d}(l_{\underline{F}}(j_1, \ldots,
     j_{\nu(\underline{F})})).
\end{eqnarray*}
    Note that the inner sum on the right-hand side of \eqref{8.1} is just $R_{{\underline{O}}}(f)$
    and observe that for any $\underline{G} \in \Pi_n$,
\begin{equation*}
    C_{{\underline{G}}}(f) = \sum_{\underline{G} \preceq
    \underline{F}}R_{{\underline{F}}}(f).
\end{equation*}
    This is merely partitioning the unrestricted sum for
    $C_{{\underline{G}}}$ as a sum over the various possibilities
    for coincidences between the indices. Thus we can use M\"obius
    inversion to express $R_{{\underline{O}}}(f)$ in terms of the
    sums $C_{{\underline{F}}}(f)$:
\begin{equation*}
    R_{{\underline{O}}}(f)=\sum_{\underline{F}}\mu(\underline{O},
    \underline{F})C_{{\underline{F}}}(f),
\end{equation*}
    so that the right-hand side of \eqref{8.1} can be expressed as
\begin{equation*}
  \frac {2^n \pi^2 }{4 X}\sum_{d \in D(X)}\sum_{\underline{F}}(-1)^{n-\nu(\underline{F})}
\left (\prod^{\nu(\underline{F})}_{l=1}(|F_l|-1)! \right
)\sum_{\substack{j_1, \cdots,j_n \\ \text{positive}}}
     \tilde{f}_{8d}\left( l_{\underline{F}}(j_1, \ldots,
     j_{\nu(\underline{F})}) \right ),
\end{equation*}
   where $\underline{F}$ ranges over all ways of decomposing $\{1, \ldots, n \}$ into disjoint subsets
$\{ F_1, $ $\ldots, F_{\nu(\underline{F})} \}$. By \eqref{01.00} and
because we are assuming that the $f_i$'s are even, we can recast the
inner sum above as going over all the $\gamma^{(j)}_{8d}$'s (
instead of $j>0$) with the presence of a factor
$1/2^{\nu(\underline{F})}$. Thus, \eqref{8.1} becomes
\begin{equation}
\label{8.1'}
  \frac {2^n \pi^2 }{4 X}\sum_{d \in D(X)}
\sum_{\underline{F}}\frac
{(-1)^{n-\nu(\underline{F})}}{2^{\nu(\underline{F})}}
\prod^{\nu(\underline{F})}_{l=1}\left( \left( |F_l|-1 \right)!
\sum_{\gamma_{8d}}\prod_{i \in F_l}f_i(L\gamma_{8d}) \right ).
\end{equation}
   We write
\begin{equation*}
   F_{l}(x)=\prod_{i \in F_l}f_i(x),
\end{equation*}
   so that we can express \eqref{8.1'} as
\begin{equation*}
  \frac {2^n \pi^2 }{4 X}\sum_{d \in D(X)}
\sum_{\underline{F}}\frac
{(-1)^{n-\nu(\underline{F})}}{2^{\nu(\underline{F})}}
\prod^{\nu(\underline{F})}_{l=1}\left( \left( |F_l|-1 \right)!
\sum_{\gamma_{8d}}F_l(L\gamma_{8d}) \right ).
\end{equation*}

   By the explicit formula \eqref{01.3}, we now arrive at the
   following (note that $\prod^{\nu(\underline{F})}_{l=1}\hat{F}_l(u_l)$ is supported in
   $\sum^{\nu(\underline{F})}_{l=1}|u_l|<2$, see Lemma \ref{lem10.1} below):
\begin{prop}
\label{prop2.5}
   Assume each $f_i$ in \eqref{01.01} is even and in $S({\mathbb
   R})$ and $\hat{f}$ is supported
   in $\sum^n_{i=1}|u_i|<2$, then
\begin{eqnarray}
\label{10.060} \label{8.3}
   &&  \frac {\pi^2 }{4 X}\sum_{d \in D(X)}{\sum_{j_1, \ldots,
     j_n}}^*f(L\gamma^{(j_1)}_{8d}, L\gamma^{(j_2)}_{8d}, \ldots,
     L\gamma^{(j_n)}_{8d})  \\
  &=& \frac {\pi^2 }{4 X}\sum_{d \in D(X)}
\sum_{\underline{F}}(-2)^{n-\nu(\underline{F})}
\prod^{\nu(\underline{F})}_{l=1}(|F_l|-1)! \left ( C_l+D_l+O\big (
\frac{\log \log X}{\log X} \big ) \right ), \nonumber
\end{eqnarray}
   where (note that $\hat{F}_l(x)$ is even since each $f_i$ is even),
\begin{eqnarray*}
   C_l &=& \int_{{\mathbb R}}F_l(x)dx, \\
   D_l &=& -\frac 1{2}\int_{{\mathbb R}}\hat{F}_l(u)du-\frac {2}{\log X}\sum_p\frac {\log
   p}{\sqrt{p}}\genfrac(){1pt}{}{8d}{p} \hat{F}_l(\frac {\log p}{\log X}).
\end{eqnarray*}
\end{prop}
   For the moment, we consider the expression in \eqref{8.3} without the term $O(\frac {\log \log X}{\log X})$.
When we expand the product over $l$, we obtain
$2^{\nu(\underline{F})}$ terms, each a product of $C_l$'s and
$D_l$'s. A typical term can be written as
\begin{equation*}
   \prod_{l \in S^c}C_l\prod_{l \in S}D_l
\end{equation*}
   for some subset $S$ of $\{1, 2, \ldots, \nu(\underline{F}) \}$. (Empty products are taken to be $1$.) The
product of the $C_l$'s contributes a factor of
\begin{equation*}
  \prod_{l \in S^c}\int_{\mathbb{R}}F_l(x)dx.
\end{equation*}

   In order to evaluate the contribution of the product of the $D_l$'s, we first
focus on finding the asymptotic expression of (with $n \geq 1$)
\begin{equation*}
\label{01}
    S(X,Y;\prod^n_{i=1}\hat{g}_i) :=
    \sum_{\substack { X \leq d \leq 2X \\ (d,2)=1}}\mu^2(d)\sum_{\prod^n_{i=1}p_i \leq Y}\left ( \prod^n_{i=1}
    \frac {\log p_i}{\sqrt {p_i}}\genfrac(){1pt}{}{8d}{p_i}\hat{g}_i(\frac {\log p_i}{\log X}) \right ),
\end{equation*}
    where $\hat{g}_i(u_i)$'s are smooth and supported on $\sum^{n}_{i=1}|u_i|< 2-\epsilon$
    with
$\epsilon$ to be specified later. To emphasis this condition, here
and throughout we shall set $Y=X^{2-\epsilon}$ and write the
condition $\prod^n_{i=1}p_i \leq Y$ explicitly. In place of $S(X,Y;
\prod^n_{i=1}\hat{g}_i)$ it is technically easier to consider the
    smoothed sum
\begin{equation*}
    S(X,Y;\prod^n_{i=1}\hat{g}_i, \Phi) :=
    \sum_{(d,2)=1}\mu^2(d)\sum_{\prod^n_{i=1}p_i \leq Y}\left ( \prod^n_{i=1}
    \frac {\log p_i}{\sqrt {p_i}}\genfrac(){1pt}{}{8d}{p_i}\hat{g}_i(\frac {\log p_i}{\log X}) \right )\Phi(\frac {d}{X}),
\end{equation*}
    where $\Phi$ is a smooth function supported on $(1,2)$,
    satisfying $\Phi(t)=1$ for $t \in (1+1/U, 2-1/U)$, and such that
    $\Phi^{(j)}(t) \ll_j U^j$ for all integers $j \geq 0$.

Let $Z >1$ be a real parameter to be chosen later and write
     $\mu^2(d)=M_Z(d)+R_Z(d)$ where
\begin{equation*}
    M_Z(d)=\sum_{\substack {l^2|d \\ l \leq Z}}\mu(l),    R_Z(d)=\sum_{\substack {l^2|d \\ l >
    Z}}\mu(l).
\end{equation*}

    Define
\begin{eqnarray*}
    S_M(X,Y;\prod^n_{i=1}\hat{g}_i,\Phi) &= &\sum_{(d,2)=1}M_Z(d)\sum_{\prod^n_{i=1}p_i \leq Y}\left ( \prod^n_{i=1}
    \frac {\log p_i}{\sqrt {p_i}}\genfrac(){1pt}{}{8d}{p_i}\hat{g}_i(\frac {\log p_i}{\log X}) \right )\Phi(\frac
    {d}{X}), \\
   S_R(X,Y;\prod^n_{i=1}\hat{g}_i, \Phi)
&=&\sum_{(d,2)=1}R_Z(d)\sum_{\prod^n_{i=1}p_i \leq Y}\left (
\prod^n_{i=1}
    \frac {\log p_i}{\sqrt {p_i}}\genfrac(){1pt}{}{8d}{p_i}\hat{g}_i(\frac {\log p_i}{\log X}) \right )\Phi(\frac {d}{X}),
\end{eqnarray*}
   so that $S(X,Y;\prod^n_{i=1}\hat{g}_i,\Phi)=S_M(X,Y;\prod^n_{i=1}\hat{g}_i,\Phi)+S_R(X,Y;\prod^n_{i=1}\hat{g}_i,\Phi)$.

    Using standard techniques, we can show that by choosing $U$ and $Z$ properly, both the terms $S_R(X,Y;\prod^n_{i=1}\hat{g}_i,
    \Phi)$ and
    $|S(X,Y;$ $\prod^n_{i=1}\hat{g}_i)$ $-$
    $S(X,Y;\prod^n_{i=1}\hat{g}_i,
    \Phi)|$  are small.
    Hence the main term arises only from
    $S_M(X,Y;\prod^n_{i=1}\hat{g}_i,\Phi)$. We write
it as
\begin{eqnarray}
\label{3.1}
   && S_{M}(X,Y;\prod^n_{i=1}\hat{g}_i, \Phi)  \\
   &=&  \sum_{\prod^n_{i=1}p_i \leq Y}\left (
\prod^k_{i=1}
    \frac {\log p_i}{\sqrt {p_i}}\genfrac(){1pt}{}{8}{p_i}\hat{g}_i(\frac {\log p_i}{\log X})
    \right
    )\left ( \sum_{(d,2)=1}M_Z(d)\genfrac(){1pt}{}{d}{\prod^n_{i=1}p_i}
    \Phi(\frac {d}{X}) \right ).  \nonumber
\end{eqnarray}

    From the above one can see that the problem of obtaining an
    asymptotic expression of the term $S_{M}(X,Y;\prod^n_{i=1}\hat{g}_i, \Phi)$
    is reduced to an evaluation of certain character sums,
    specifically,
\begin{equation*}
    \sum_{(d,2)=1}M_Z(d)\genfrac(){1pt}{}{d}{\prod^n_{i=1}p_i}
    \Phi(\frac {d}{X}).
\end{equation*}

    To help understanding the situation, we consider a simple case,
    $Z=1$ here and we also replace $\prod^n_{i=1}p_i$ with $P$ in the above
    sum, with no constraint on the value of $P$. Hence we are lead
    to consider the smoothed character sum
\begin{equation}
\label{2.22}
    \sum_{(d,2)=1}\genfrac(){1pt}{}{8d}{P}
    \Phi(\frac {d}{X}).
\end{equation}

     We certainly expect the above sum is large when $P=\Box$, in
     which case the character $\genfrac(){1pt}{}{\cdot}{P}$ is principal.

     When $P \neq \Box$, we note that the length of our character
     sum is of size $X$ and when $P \leq X$, we are dealing with a
     long character sum and therefore we expect substantial cancellations. This is exactly what was carried out
     in Rubinstein's work \cite{R}.

     When $P \neq \Box$ and $X < P$, we are dealing with a
     short character sum.
     Our idea then is to apply the Poisson summation formula (see Lemma \ref{lem2} below) to convert
     \eqref{2.22} to another sum:
\begin{equation*}
   \sum_{(d,2)=1}\genfrac(){1pt}{}{8d}{P}
    \Phi(\frac {d}{X})=\frac {X}{2P}\genfrac(){1pt}{}{16}{P}
       \sum_m(-1)^mG_m(P)\tilde{\Phi}(\frac {mX}{2P}),
\end{equation*}
    with the definition of $G_m$ and $\tilde{\Phi}$ given in
    Lemma \ref{lem1} and \eqref{012}, respectively. Here
    $\tilde{\Phi}(x)$ can be thought as being supported in $|x|<1$
    and $G_m(P)$ can be treated as a character. Hence we have
    converted \eqref{2.22} to another character sum with length of size $P/X$
    and we expect to obtain a better estimation this way as long as $P/X
    < X$ or $P < X^2$, which is exactly why we can have $\sum^n_{i=1}|u_i|<2$ in Theorem
     \ref{thm2.1}.

\section{The Proof of Theorem \ref{thm2.1}}
\label{sec3.0} \setcounter{equation}{0}

\subsection{The Term $S_R$}
\label{sec 3.1} \setcounter{equation}{0}
   We recall here in Chapter \ref{sec2.1} we have defined
\begin{equation*}
   S_R(X,Y;\prod^n_{i=1}\hat{g}_i, \Phi)
= \sum_{(d,2)=1}R_Z(d)\sum_{\prod^n_{i=1}p_i \leq Y}\left (
\prod^n_{i=1}
    \frac {\log p_i}{\sqrt {p_i}}\genfrac(){1pt}{}{8d}{p_i}\hat{g}_i(\frac {\log p_i}{\log X}) \right )\Phi(\frac
    {d}{X}).
\end{equation*}
   We will show in this section that this is negligible. As a preparation,
   we first seek a bound for
\begin{equation*}
  E(U;\chi,\prod^n_{i=1}\hat{g}_i)=\sum_{\prod^n_{i=1}p_i \leq U}\left (
\prod^n_{i=1}
    \frac {\log p_i}{\sqrt {p_i}}\chi (p_i)\hat{g}_i(\frac {\log p_i}{\log X})
    \right
    ),
\end{equation*}
   for any non-principal quadratic character $\chi$ with modulus $q$ and $U \leq
   X^2$. We have the following
\begin{prop}
\label{prop3.1}
   Assume GRH. For any non-principal quadratic character $\chi$ with modulus $q$ and $U \leq
   X^2$ and with $\prod^n_{i=1}\hat{g}_i$ smooth and supported in $\sum^n_{i=1}|u_i| <2$,
\begin{equation*}
  E(U;\chi,\prod^n_{i=1}\hat{g}_i) \ll \log^{2n+1} (qX).
\end{equation*}
\end{prop}
\begin{proof}
  Both the estimation and the method we use for the proof are
  standard and we include a proof here for completeness.
  We first assume $2 \leq U \leq X^2$ and note the following discontinuous integral
   \cite[\S 17]{D}:
\begin{equation*}
    \frac 1{2 \pi i} \int_{(c)}\frac {y^s}{s}ds=\left\{\begin{array}{cl}
    0  & \mbox{if $ 0 < y <1 $}, \\
    \frac {1}{2} & \mbox{if $y=1$},  \\
    1  & \mbox{if $y>1$},
    \end{array}\right.
\end{equation*}
    where $c>0$ and that \cite[\S 17]{D} for $y>0, c>0, T>0$,
\begin{equation*}
    \int_{(c)}\frac {y^s}{s}ds- \int^{c+iT}_{c-iT}\frac
    {y^s}{s}ds \ll \left\{\begin{array}{cl}
    y^c\min (1, T^{-1}|\log y|^{-1} ) & \mbox{if $ y \neq 1 $}, \\
    cT^{-1} & \mbox{if $y=1$}.
    \end{array}\right.
\end{equation*}

    Using the above, we can recast $E(U;\chi, \prod^n_{i=1}\hat{g}_i)$ as
\begin{eqnarray}
\label{3.01}
  &&    E(U;\chi,\prod^n_{i=1}\hat{g}_i)  \\
 &=&\frac {1}{2 \pi i} \int^{c+iT}_{c-iT} \prod^n_{i=1}
\left (\sum_{p}
    \frac {\log p}{p^{s+1/2}}\chi (p)\hat{g}_i(\frac {\log p}{\log X})
    \right
    )U^s \frac {ds}{s} \nonumber \\
    && + O\Big( (\frac c{T}+1) \prod^n_{\substack{i=1 \\ \prod^n_{i=1}p_i=U}}
    \frac {\log p_i}{\sqrt{p_i}} \Big ) \nonumber \\
    && + O \left (\sum_{\prod^n_{i=1}p_i \leq X^2}\left (
\prod^n_{i=1}
    \frac {\log p_i}{\sqrt {p_i}} \right
    )(\frac {U}{\prod^n_{i=1}p_i})^c\min (1, T^{-1}|\log \frac {U}{\prod^n_{i=1}p_i}|^{-1} ) \right
    ). \nonumber
\end{eqnarray}

    For the first error term above, we may assume $U$ is an integer and we use $\log p/\sqrt{p} \ll
    p^{-1/2+\delta}$ for an arbitrarily small number $\delta$ so that the each single product is of size
    $\ll U^{-1/2+\delta}$ and there are at most $(\omega(U))^n$
    choices for selecting the $p_i$'s. Here $\omega(m)$ denotes the
    number of distinct primes dividing $m$ and it's well-known that
    for $m \geq 3$,
\begin{equation*}
    \omega(m) \ll \frac {\log m}{ \log \log m}.
\end{equation*}
     It follows from this
    that the error term is of size
\begin{equation*}
   \ll (\frac c{T}+1) U^{-1/2+\delta}.
\end{equation*}

    Note that $|\log (U/\prod^n_{i=1}p_i)|$ is bounded below for $\prod^n_{i=1}p_i
    \geq \frac {5}{4}U$ or $\prod^n_{i=1}p_i
    \leq \frac {3}{4}U$ and those $\prod^n_{i=1}p_i$ in this range contribute
    to the second error term above
\begin{equation*}
    \ll \frac {U^cX \log^{n-1}X}{T}.
\end{equation*}

    For the terms $\frac {3}{4}U < \prod^n_{i=1}p_i < U$, we note that
\begin{equation*}
    \log \frac {U}{\prod^n_{i=1}p_i}=-\log (1-\frac {U-\prod^n_{i=1}p_i}{U}) \geq \frac
    {U-\prod^n_{i=1}p_i}{U}.
\end{equation*}
   Hence these terms contribute
\begin{equation*}
    \left(\max_{3U/4 \leq m \leq U}\omega(m) \right)^n(\frac {4}{3})^c\frac {U}{T}\sum^{U/4}_{i=1}\frac 1{i}
    \ll (\frac {4}{3})^c\frac {U \log^{n+1}
    U}{T}.
\end{equation*}

   Similarly, the terms $U < \prod^n_{i=1}p_i < \frac {5}{4}U$ contribute
\begin{equation*}
    \ll \frac {U \log^{n+1}  U}{T}.
\end{equation*}

   It remains to treat the integral on the right-hand side of \eqref{3.01}.
   We first look at the case when $\chi$ is primitive. Note for any fixed $|t| \leq T$, $V \geq 2$,
\begin{equation*}
    \sum_{p \leq V}\frac {\log p}{p^{it}}\chi(p)=\sum_{m
\leq V}\frac {\Lambda(m)}{m^{it}}\chi(m)+O(V^{1/2}\log V).
\end{equation*}
    Now write for $c_1=1+1/\log V$,
\begin{eqnarray*}
    \sum_{m
\leq V}\frac {\Lambda(m)}{m^{it}}\chi(m) &=& \frac 1{2 \pi
i}\int_{(c_1)}\sum^{\infty}_{m=1}\frac
{\Lambda(m)}{m^{s+it}}\chi(m)V^s \frac {ds}{s}+O(\log V) \\
  &=& \frac 1{2 \pi
i}\int_{(c_1)}-\frac {L'(s+it, \chi)}{L(s+it, \chi)}V^s \frac
{ds}{s}+O(\log V).
\end{eqnarray*}
   Similar as in \cite[\S 19]{D} and assume GRH, we have for $T_1>|t|+2$,
\begin{eqnarray*}
  &&  \frac 1{2 \pi
i}\int_{(c_1)}-\frac {L'(s+it, \chi)}{L(s+it, \chi)}V^s \frac
{ds}{s} \\
&=&\frac 1{2 \pi i}\int^{c_1+iT_1}_{c_1-iT_1}-\frac {L'(s+it,
\chi)}{L(s+it, \chi)}V^s \frac {ds}{s}+O(\frac {V \log^2
V}{T_1}+\log V) \\
  &=& -\sum_{|\gamma-t| < T_1}\frac {V^{\rho-it}}{\rho-it} +O(\frac {V\log^2 (q(T_1-|t|))}{(T_1-|t|)\log V}+\frac {V \log^2
V}{T_1}+\log V+\log^2 q )
\end{eqnarray*}
  where $\rho=1/2+i\gamma$ denotes a non-trivial zero of the
$L$-function. Let $N(T,\chi)$ denote the number of zeros of $L(s,
\chi)$ in the rectangle $0<\sigma<1, |t|<T$. By \cite[\S 16]{D},
$N(T+1,\chi)-N(T, \chi) \ll \log (qT)$ so that
\begin{equation*}
   \sum_{|\gamma-t| < T_1}\frac {V^{\rho-it}}{\rho-it} \ll
V^{1/2}\log \left(q \left(T_1+|t| \right )\right)\sum_{i \leq
2T_1}\frac 1{i} \leq V^{1/2}\log^2 \left(q \left(T_1+|t| \right
)\right).
\end{equation*}

   For our application here, we will take $V \leq T=X^2$. Hence upon taking $T_1=X^4$, we obtain
\begin{equation*}
\sum_{p \leq V}\frac {\log p}{p^{it}}\chi(p) \ll V^{1/2}\log^2(qX).
\end{equation*}
   The above holds also for $1 \leq V < 2$, hence by partial summation, for $\chi$ primitive,
   $c=1/\log X,|t| \leq T, 1 \leq i \leq n$,
\begin{equation}
\label{3.02}
  \sum_{p_i}\frac {\log
p_i}{p^{1/2+c+it}_i}\chi(p)\hat{g}_i(\frac {\log p_i}{\log X}) \ll
\log^2(qX).
\end{equation}

   Now for $\chi$ induced by $\chi_1$ with modulus
$q_1|q$, we have
\begin{eqnarray*}
  && \sum_{p_i}\frac {\log
p_i}{p^{1/2+c+it}_i}\chi(p)\hat{g}_i(\frac {\log p_i}{\log
X})-\sum_{p_i}\frac {\log
p_i}{p^{1/2+c+it}_i}\chi_1(p)\hat{g}_i(\frac {\log p_i}{\log X}) \\
&\ll& \sum_{p_i | q}\log p_i \ll \log q.
\end{eqnarray*}
  Hence \eqref{3.02} still holds in this case and on setting $c=1/\log X, T=X^2$, we obtain for
   $2 \leq U \leq X^2$,
\begin{equation*}
    E(U;\chi,\prod^n_{i=1}\hat{g}_i)  \ll  \log^{2n} (qX)  \int^{T}_{1} \frac {dt}{t}+O(\log^{2n}(qX)\log X)
\ll \log^{2n+1} (qX).
\end{equation*}
    Note that $E(U;\chi,\prod^n_{i=1}\hat{g}_i)=0$ when $1
 \leq U <2$ and this completes the proof.
\end{proof}

   As a consequence of Proposition \ref{prop3.1}, we have
\begin{lemma}
\label{lem3.1}
   With the definition at the beginning of this section,
\begin{equation*}
    S_R(X,Y;\prod^n_{i=1}\hat{g}_i, \Phi) \ll \frac {X \log^{2n+1} X}{Z}.
\end{equation*}
\end{lemma}
\begin{proof}
  On writing $d=l^2m $, we obtain
\begin{eqnarray*}
    S_R(X,Y;\prod^n_{i=1}\hat{g}_i, \Phi) &=&
     \sum_{\substack{l>Z \\ (l, 2)=1}}\mu(l)\sum_{(m,2)=1}
     \Phi(\frac
     {l^2m}{X})E(Y;\chi_{8l^2m},\prod^n_{i=1}\hat{g}_i)
     \\
     &\ll&  \sum_{l>Z}\sum_{X/l^2 \leq m \leq
   2X/l^2} \log^{2n+1}(X) \ll \frac {X \log^{2n+1} X}{Z}.
\end{eqnarray*}
\end{proof}

\subsection{The Term $S_M$}
\label{sec 3.2}
   As explained in Section \ref{sec 2.3},
   we will evaluate $S_{M}(X,Y;\prod^n_{i=1}{\hat g}_i, \Phi)$ by applying the Poisson summation
    formula to the sum over $d$ in $S_{M}(X,Y;\prod^n_{i=1}{\hat g}_i,\Phi)$. For all odd
    integers $k$ and all integers $m$, we introduce the Gauss-type
    sums
\begin{equation*}
\label{010}
    \tau_m(k) := \sum_{a ~~(\mbox{mod}~~ k)}(\frac {a}{k})e(\frac {am}{k}) =:
    (\frac {1+i}{2}+(\frac {-1}{k})\frac {1-i}{2})G_m(k),
\end{equation*}
   where $e(x):=e^{2\pi i x}$ as defined before. We quote Lemma 2.3 of
   \cite{sound1} which determines $G_m(k)$.
\begin{lemma}
\label{lem1}
   If $(k_1,k_2)=1$ then $G_m(k_1k_2)=G_m(k_1)G_m(k_2)$. Suppose that $p^a$ is
   the largest power of $p$ dividing $m$ (put $a=\infty$ if $m=0$).
   Then for $b \geq 1$ we have
\begin{equation*}
\label{011}
    G_m(p^b)= \left\{\begin{array}{cl}
    0  & \mbox{if $b\leq a$ is odd}, \\
    \phi(p^b) & \mbox{if $b\leq a$ is even},  \\
    -p^a  & \mbox{if $b=a+1$ is even}, \\
    (\frac {m/p^a}{p})p^a\sqrt{p}  & \mbox{if $b=a+1$ is odd}, \\
    0  & \mbox{if $b \geq a+2$}.
    \end{array}\right.
\end{equation*}
\end{lemma}

   For $F \in S(\mathbb{R})$ we define
\begin{equation}
\label{012}
   \tilde{F}(\xi)=\frac {1+i}{2}\hat{F}(\xi)+\frac
   {1-i}{2}\hat{F}(-\xi)=\int^{\infty}_{-\infty}\left(\cos(2\pi \xi
   x)+\sin(2\pi \xi x) \right)F(x)dx.
\end{equation}
    We quote Lemma 2.6 of \cite{sound1} which determines the
   inner sum in \eqref{3.1}.
\begin{lemma}
\label{lem2}
   Let $\Phi$ be a non-negative, smooth function supported in
   $(1,2)$. For any odd integer $k$,
\begin{equation*}
\label{013}
  \sum_{(d,2)=1}M_Z(d)(\frac {d}{k})
    \Phi(\frac {d}{X})=\frac {X}{2k}(\frac {2}{k})\sum_{\substack {\alpha \leq Z \\ (\alpha, 2k)=1}}\frac {\mu(\alpha)}{\alpha^2}
    \sum_m(-1)^mG_m(k)\tilde{\Phi}(\frac {mX}{2\alpha^2k}).
\end{equation*}
\end{lemma}

   Now we can write $S_{M}(X,Y;\prod^n_{i=1}{\hat g}_i, \Phi)$ as
\begin{eqnarray*}
\label{014}
  S_{M}(X,Y;\prod^n_{i=1}{\hat g}_i, \Phi) &=& \frac {X}{2}\sum_{\substack{\prod^n_{i=1}p_i \leq Y \\ (2,
\prod^n_{i=1}p_i)=1
 }}\left ( \prod^{n}_{i=1}\frac {\log p_i}{p^{3/2}_i}{\hat g}_i(\frac {\log p_i}{\log
 X}) \right ) \\
 && \cdot \sum_{\substack {\alpha \leq Z \\ (\alpha,
 2\prod^n_{i=1}p_i)=1 }} \frac {\mu(\alpha)}{\alpha^2}
    \sum_m(-1)^mG_m(\prod^n_{i=1}p_i)\tilde{\Phi}(\frac {m X}{2\alpha^2\prod^n_{i=1}p_i}).
\end{eqnarray*}

\subsection{The First Main Term}
\label{sec 3.3}
    We consider the sum $S_{0}$ in $S_{M}(X,Y;\prod^n_{i=1}{\hat g}_i, \Phi)$ corresponding to the contribution of $m=0$.
It follows directly from the
    definition that $G_0(k)=\varphi(k)$ if $k=\Box$ and $G_0(k)=0$
    otherwise. Thus $S_{0}=0$ when $n$ is odd and for the even
    $n$'s,
\begin{equation*}
    S_{0} = \frac {X\hat{\Phi}(0)}{2}\sum_{\substack{\prod^n_{i=1}p_i \leq Y \\ (2, \prod^n_{i=1}p_i)=1
 \\ \prod^n_{i=1}p_i=\Box}}\left ( \prod^{n}_{i=1}\frac {\log p_i}{p^{3/2}_i}{\hat g}_i(\frac {\log p_i}{\log
 X}) \right )\phi(\prod^n_{i=1}p_i)\sum_{\substack {\alpha \leq Z \\ (\alpha,
 2\prod^n_{i=1}p_i)=1}}\frac {\mu(\alpha)}{\alpha^2}.
\end{equation*}
    For a partition of the set $\{1, 2, \ldots, n \}$ into $k$
    subsets $\{ I_i \}$ with $|I_i|=\alpha_i$, we denote
    $$ h_i(\frac {\log p}{\log X})=\prod_{j \in I_i}{\hat g}_j(\frac {\log p}{\log X}).$$
   Using the above notation, we can express $S_0$ as linear combinations of the
    sums of the following form
\begin{eqnarray}
\label{3.002}
  && \sum_{\substack {\prod^{k}_{i=1}p^{\alpha_i}_i \leq Y \\ (2, \prod^k_{i=1}p_i)=1 \\  p_i \neq p_j \\ 2 | \alpha_i}}
   \prod^{k}_{i=1}\frac {(\log p_i)^{\alpha_i}}{p^{3\alpha_i/2}_i}
   h_i(\frac {\log p_i}{\log X})\phi(p^{\alpha_i}_i) \sum_{\substack { \alpha \leq Z \\ (\alpha,
2\prod^k_{i=1}p_i)=1}}\frac {\mu(\alpha)}{\alpha^2} \\
 &=& \frac 8{\pi^2}\sum_{\substack {\prod^{k}_{i=1}p^{\alpha_i}_i \leq Y \\
(2, \prod^k_{i=1}p_i)=1 \\  p_i \neq p_j \\ 2 | \alpha_i}}
   \prod^{k}_{i=1}\frac {(\log p_i)^{\alpha_i}}{p^{\alpha_i/2}_i}
   h_i(\frac {\log p_i}{\log X})(1+\frac 1{p_i})^{-1}\left(1+O(\frac 1{Z})
\right ),  \nonumber
\end{eqnarray}
  where we have used
\begin{equation*}
    \sum_{\substack{ \alpha \leq Z \\ (\alpha,
m)=1}}\frac {\mu(\alpha)}{\alpha^2}=\frac
1{\zeta(2)}\prod_{p|m}(1-\frac 1{p^2})^{-1}\left(1+O(\frac 1{Z})
\right ).
\end{equation*}
   If $\alpha_i \geq 4$ for some $i$ then the right-hand side
   expression in \eqref{3.002} is readily seen to be
\begin{equation*}
  \ll \log^{n-4} X.
\end{equation*}
   Since there are $O_n(1)$ partitions of the set $\{ 1, \ldots, n
   \}$, we can now write $S_0$ as
\begin{equation*}
  S_0 = \frac {4X\hat{\Phi}(0)}{\pi^2} \sum_{\substack{(A;B) \\ \prod^{n/2}_{i=1}p^2_{i} \leq Y \\ p_i \neq p_j}}
\left ( \prod^{n/2}_{i=1}\frac {\log^2 p_i}{p_i}{\hat g}_{a_i}(\frac
{\log p_i}{\log
 X}){\hat g}_{b_i}(\frac {\log p_i}{\log
 X})(1+\frac 1{p_i})^{-1} \right )\left(1+O(\frac 1{Z})
\right ),
\end{equation*}
   where $\sum_{(A;B)}$ is over all ways of pairing up the elements of $\{1,\ldots,
   n\}$. By using similar arguments as above, we see that removing
   the condition $p_i \neq p_j$ introduces an error term of size
   $\ll X \log^{n-1}X$, hence
\begin{eqnarray}
\label{4.1}
    S_{0} &=&
      \frac {4X\hat{\Phi}(0)}{\pi^2} \sum_{\substack{(A;B) \\ \prod^{n/2}_{i=1}p^2_{i} \leq Y}}
\left ( \prod^{n/2}_{i=1}\frac {\log^2 p_i}{p_i}{\hat g}_{a_i}(\frac
{\log p_i}{\log
 X}){\hat g}_{b_i}(\frac {\log p_i}{\log
 X})(1+\frac 1{p_i})^{-1}  \right )   \\
 && \cdot \left(1+O(\frac 1{Z})
\right )+O( X \log^{n-1}X) \nonumber \\
  &=&  \frac {4X\log^n X \hat{\Phi}(0)}{\pi^2}\sum_{(A;B)}\prod^{n/2}_{i=1}
\int^{\infty}_0 u_i{\hat g}_{a_i}(u_i){\hat g}_{b_i}(u_i)du_i  \nonumber \\
 &&  +O( X \log^{n-1}X+ \frac {X\log^n X}{Z}), \nonumber
\end{eqnarray}
   where the integral on the right-hand side of \eqref{4.1} comes from
partial summation and Mertens' formula \cite[p. 57]{D}
\begin{equation*}
    \sum_{p \leq x}\frac {\log p}{p}= \log x+O(1),
\end{equation*}
.

\subsection{The Error Term}
\label{sec 3.4}
    We now consider the sums in $S_{M}(X,Y;\prod^n_{i=1}{\hat g}_i, \Phi)$ corresponding to the
    contribution of $m \neq 0, \Box$. Using the notations in Section \ref{sec
    3.3}, the sums in $S_{M}(X,Y;\prod^n_{i=1}{\hat g}_i,
\Phi)$ corresponding to the
    contribution of $m \neq 0, \Box$ are linear combinations of the
    sums of the following form
\begin{eqnarray*}
  && \sum_{\substack {\prod^{k}_{i=1}p^{\alpha_i}_i \leq Y \\ (2, \prod^k_{i=1}p_i)=1 \\  p_i \neq p_j}}
   \prod^{k}_{i=1}\frac {(\log p_i)^{\alpha_i}}{p^{3\alpha_i/2}_i}
   h_i(\frac {\log p_i}{\log X}) \\
 && \cdot \sum_{\substack { \alpha \leq Z \\ (\alpha,
2\prod^k_{i=1}p_i)=1}}\frac {\mu(\alpha)}{\alpha^2}
    \sum_{m \neq 0, \Box}(-1)^mG_m(p^{\alpha_i}_i) \tilde{\Phi}
    (\frac {m X}{2\alpha^2\prod^{k}_{i=1}p^{\alpha_i}_i}).
\end{eqnarray*}
    There are $O_n(1)$ partitions of the set $\{ 1, \ldots, n \}$ and we now consider the above sum for
   a fixed partition. Without loss of generality, we may assume
   there exists integers $k_1, k_2$ with $0 \leq k_2 \leq k_1 \leq k$
  such that
   for
   $k_2 < i \leq k$, $\alpha_i \equiv 1 \hspace{0.05in}(\text{mod } 2)$
    and $\alpha_i$ $\equiv 0 \hspace{0.05in} (\text{mod } 2)$ for $i \leq
   k_2$; for $k_1 < i \leq k$, $\alpha_i = 1$ and for $k_2 < i \leq k_1$, $\alpha_i \geq
   3$.
   By Lemma \ref{lem1},
   each $m$ can be written as $m=\prod^{k}_{i=k_2+1}p^{\alpha_i-1}_im'$ and we take empty product to be $1$.
   Hence with a slightly change of notation, the sum we are interested
   in can be written as
\begin{eqnarray*}
  R &=& \sum_{\substack { \alpha \leq Z \\ (\alpha,
2)=1}}\frac {\mu(\alpha)}{\alpha^2}\sum_{\substack
{\prod^{k}_{i=1}p^{\alpha_i}_i \leq Y \\ (2\alpha, \prod^k_{i=1}p_i)=1 \\
p_i \neq p_j}}
   \prod^{k_2}_{i=1}\frac {(\log p_i)^{\alpha_i}}{p^{3\alpha_i/2}_i}
   h_i(\frac {\log p_i}{\log X}) \nonumber
   \\
   && \cdot \prod^{k_1}_{i=k_2+1}\frac {(\log p_i)^{\alpha_i}}{p^{(\alpha_i+1)/2}_i}
  h_i(\frac {\log p_i}{\log X})
  \prod^{k}_{i=k_1+1}(\frac {\log p_i}{p_i})
  h_i(\frac {\log p_i}{\log X}) \nonumber \\
  &&   \cdot \sum_{m \neq 0, \Box}(-1)^m\tilde{\Phi}
    (\frac {m X}{2\alpha^2\prod^{k_2}_{i=1}p^{\alpha_i}_i \prod^{k}_{i=k_2+1}p_i})
    \genfrac(){1pt}{}{m}{\prod^{k}_{k_2+1}p_i}\Big (\prod^{k_2}_{i=1}G_m(p^{\alpha_i}_i) \Big ). 
\end{eqnarray*}
   For any non-principal quadratic character $\chi$ with modulus $q$
   and for $\prod^n_{i=1}\hat{g}_i$ smooth with support in $\sum^n_{i=1}|u_i|
<2$, we write for $U \leq X^2$,
\begin{equation*}
    D(U;\chi,\prod^m_{i=1}{\hat g}_i)=\sum_{\substack { \prod^m_{i=1}p_i \leq U \\ p_i \neq p_j}}\left (
\prod^m_{i=1}
    \frac {\log p_i}{\sqrt {p_i}}\chi (p_i){\hat g}_i(\frac {\log p_i}{\log X})
    \right
    ),
\end{equation*}
   and similar to our treatment of $E(U;\chi,\prod^n_{i=1}{\hat g}_i)$ in
   Section \ref{sec 3.1}, one has
\begin{equation}
\label{3.04}
    D(U;\chi,\prod^m_{i=1}{\hat g}_i) \ll \log^{2m+1}(qX).
\end{equation}
    Note that
\begin{equation*}
    \genfrac(){1pt}{}{4\alpha^2(\prod^{k_1}_{i=1}p_i)^2}{\prod^k_{i=k_1+1}p_i}
    = \left\{\begin{array}{cl}
    0  & \mbox{if $(2\alpha\prod^{k_1}_{i=k1}p_i, \prod^k_{i=k_1+1}p_i)>1$}, \\
    1 & \mbox{if $(2\alpha\prod^{k_1}_{i=1}p_i, \prod^k_{i=k_1+1}p_i)=1$}.
    \end{array}\right.
\end{equation*}
   so that we can express the condition $(2\alpha\prod^{k_1}_{i=1}p_i, \prod^k_{i=k_1+1}p_i)=1$
   by using the the character
   $\chi_{4\alpha^2(\prod^{k_1}_{i=1}p_i)^2}$. As
   $\chi_{4\alpha^2(\prod^{k_1}_{i=1}p_i)^2m}$ is a non-principal quadratic
   character and note that $\prod^{k}_{i=1}h(u_i)$ is supported in
   $\sum^{k}_{i=1}|u_i|<2$ (see Lemma \ref{lem10.1} below), we obtain by using \eqref{3.04} and partial summation,
   that
\begin{eqnarray*}
 && \sum_{\substack
{\prod^{k}_{i=k_1+1}p_i \leq Y/\prod^{k_1}_{i=1}p^{\alpha_i}_i \\ (2\alpha\prod^{k_1}_{i=1}p_i, \prod^k_{i=k_1+1}p_i)=1 \\
p_i \neq p_j}} \left ( \prod^{k}_{i=k_1+1}(\frac {\log p_i}{p_i})
  h_i(\frac {\log p_i}{\log X}) \right )  \\
  && \cdot \tilde{\Phi}
    (\frac {m X}{2\alpha^2\prod^{k_2}_{i=1}p^{\alpha_i}_i \prod^{k}_{i=k_2+1}p_i})
    \genfrac(){1pt}{}{m}{\prod^{k}_{k_1+1}p_i} \\
  &=& \int^{Y/\prod^{k_1}_{i=1}p^{\alpha_i}_i}_1\frac {1}{\sqrt{V}}\tilde{\Phi}
    (\frac {m X}{2\alpha^2\prod^{k_2}_{i=1}p^{\alpha_i}_i
    \prod^{k_1}_{i=k_2+1}p_iV})dD(V;\chi_{4\alpha^2(\prod^{k_1}_{i=1}p_i)^2m},\prod^k_{i=k_1+1}h_i)
    \\
   & \ll & \log^{2n+1}(X)\log^{2n+1}(|m|+2)
    \Biggr ( \frac {\sqrt{\prod^{k_1}_{i=1}p^{\alpha_i}_i}}{\sqrt{Y}}
   \left |\tilde{\Phi}\left (\frac {m\prod^{k_1}_{i=k_2+1}p^{\alpha_i-1}_i
   X}{2\alpha^2Y} \right ) \right |   \\
   && +  \int^{Y/\prod^{k_1}_{i=1}p^{\alpha_i}_i}_{1}\frac
   {1}{V^{3/2}}\left |\tilde{\Phi}
    \left (\frac {m X}{2\alpha^2\prod^{k_2}_{i=1}p^{\alpha_i}_i
    \prod^{k_1}_{i=k_2+1}p_iV} \right ) \right |dV \\
   &&   + \int^{Y/\prod^{k_1}_{i=1}p^{\alpha_i}_i}_{1}\frac
   {X}{\alpha^2(\prod^{k_2}_{i=1}p^{\alpha_i}_i
    \prod^{k_1}_{i=k_2+1}p_i) V^{5/2}} \\
   && \cdot \left |m\tilde{\Phi}'
    \left(\frac {m X}{2\alpha^2\prod^{k_2}_{i=1}p^{\alpha_i}_i
    \prod^{k_1}_{i=k_2+1}p_iV} \right ) \right |dV \Biggr ).
\end{eqnarray*}

    Hence we have
\begin{equation*}
  R  \ll \sum_{\alpha \leq Z}\frac {\log^{2n+1}(X)}{\alpha^2}(R_1+R_2+R_3),
\end{equation*}
    where
\begin{eqnarray*}
  R_1 & = & \frac {1}{\sqrt{Y}}\sum_{
  \prod^{k_1}_{i=1}p^{\alpha_i}_i \leq Y}
   \prod^{k_2}_{i=1}\frac {(\log p_i)^{\alpha_i}}{p^{\alpha_i}_i}
      \prod^{k_1}_{i=k_2+1}\frac {(\log p_i)^{\alpha_i}}{p^{1/2}_i} \\
   && \cdot \sum_{m \neq 0, \Box}\log^{2n+1}(|m|+2)\left |\tilde{\Phi}\left (\frac {m\prod^{k_1}_{i=k_2+1}p^{\alpha_i-1}_i
   X}{2\alpha^2Y}\right ) \right | \Big (\prod^{k_2}_{i=1}|G_m(p^{\alpha_i}_i)| \Big ), \\
   R_2 &=& \int^{Y}_{1}\frac
   {1}{V^{3/2}}\sum_{
  \prod^{k_1}_{i=1}p^{\alpha_i}_i \leq Y}
   \prod^{k_2}_{i=1}\frac {(\log p_i)^{\alpha_i}}{p^{3\alpha_i/2}_i}
   \prod^{k_1}_{i=k_2+1}\frac {(\log p_i)^{\alpha_i}}{p^{(\alpha_i+1)/2}_i} \\
   && \cdot \sum_{m \neq 0, \Box}\log^{2n+1}(|m|+2) \left |\tilde{\Phi}
    \left (\frac {m X}{2\alpha^2\prod^{k_2}_{i=1}p^{\alpha_i}_i
    \prod^{k_1}_{i=k_2+1}p_iV} \right ) \right |\Big (\prod^{k_2}_{i=1}|G_m(p^{\alpha_i}_i)| \Big )dV  \\
   R_3 &=& \int^{Y}_{1}\frac
   {X}{\alpha^2V^{5/2}}\sum_{
  \prod^{k_1}_{i=1}p^{\alpha_i}_i \leq Y}
   \prod^{k_2}_{i=1}\frac {(\log p_i)^{\alpha_i}}{p^{5\alpha_i/2}_i}
     \prod^{k_1}_{i=k_2+1}\frac {(\log p_i)^{\alpha_i}}{p^{(\alpha_i+3)/2}_i}  \\
   && \cdot \sum_{m \neq 0, \Box}\log^{2n+1}(|m|+2) \left|m\tilde{\Phi}'
    \left(\frac {m X}{2\alpha^2\prod^{k_2}_{i=1}p^{\alpha_i}_i
    \prod^{k_1}_{i=k_2+1}p_iV} \right) \right |\Big (\prod^{k_2}_{i=1}|G_m(p^{\alpha_i}_i)| \Big )dV.
\end{eqnarray*}

    Now we need a lemma
\begin{lemma}
\label{lemma3.5}
   For $N>0, M \geq 2$,
\begin{eqnarray}
\label{3.05}
   \sum_{m \neq 0, \Box} \log^{2n+1}(|m|M)\left |\tilde{\Phi}(\frac
   {mM}{N}) \right | &\ll& \Bigl (\log^{2n+1}M(N+2) \Bigr ) \frac
   {UN}{M}, \\
\label{3.06}
   \sum_{m \neq 0, \Box} \log^{2n+1}(|m|M) \left |mM\tilde{\Phi}'(\frac
   {mM}{N}) \right | &\ll& \Bigl(\log^{2n+1}M(N+2) \Bigr ) \frac
   {U^2N^2}{M}.
\end{eqnarray}
\end{lemma}
\begin{proof}
    Observe for $|\xi| < 1$,
\begin{equation*}
\label{16}
    \tilde{\Phi}(\xi)=\int^{2}_{1}\Bigl (\cos \left (2\pi \xi
   x)+\sin(2\pi \xi x \right ) \Bigr )dx+O(\frac 1{U}) \ll 1+O(\frac 1{U}),
\end{equation*}
    and similarly for any $i \geq 0$,
\begin{equation}
\label{16'}
    \tilde{\Phi}^{(i)}(\xi)\ll 1, \hspace{0.1in} |\xi| < 1.
\end{equation}
     Also note via integration by parts,
\begin{equation*}
\label{17}
    \tilde{\Phi}(\xi)=\frac {-1}{2\pi \xi}\Bigl(\int^{1+1/U}_{1}+\int^{2-1/U}_{2} \Bigr )
    \tilde{\Phi}'(\xi) \Bigl(\sin(2\pi \xi
   x)-\cos(2\pi \xi x) \Bigr )dx \ll \frac 1{|\xi|}.
\end{equation*}
    Similarly, one can show for any $i \geq 0, j \geq 1$,
\begin{equation}
\label{18}
    \tilde{\Phi}^{(i)}(\xi) \ll \frac {U^{j-1}}{|\xi|^j}.
\end{equation}
   Now using \eqref{16'} with $i=0$ and \eqref{18} with $i=0, j=2$ we
obtain
\begin{eqnarray*}
   &&\sum_{m \neq 0, \Box} \log^{2n+1}(|m|M) \left |\tilde{\Phi}(\frac
   {mM}{N}) \right |  \\
   &\ll& \sum_{0< |mM| \leq N} \Bigl(\log^{2n+1}M+ \log^{2n+1}(N+2) \Bigr ) \\
   && + \sum_{|mM| > N} \Bigl(\log^{2n+1}M+
   \log^{2n+1}(|m|+2)\Bigl)\frac {UN^2}{m^2M^2} \\
   &\ll & \Bigl(\log^{2n+1}M+\log^{2n+1}(N+2) \Bigr ) \frac
   {UN}{M},
\end{eqnarray*}
   and this proves \eqref{3.05}. Similarly, using \eqref{16'} with $i=0$ and \eqref{18} with $i=1, j=3$ we
obtain
\begin{eqnarray*}
   &&\sum_{m \neq 0, \Box} \log^{2n+1}(|m|M) \left|mM\tilde{\Phi}'(\frac
   {mM}{N}) \right |  \\
   &\ll& \sum_{0< |mM| \leq N} \Bigl(\log^{2n+1}M+ \log^{2n+1}(N+2) \Bigr )|mM| \\
   && + \sum_{|mM| > N} \Bigl(\log^{2n+1}M+
   \log^{2n+1}(|m|+2)\Bigl)\frac {U^2N^3}{m^2M^2} \\
   &\ll & \Bigl(\log^{2n+1}M+\log^{2n+1}(N+2) \Bigr ) \frac
   {U^2N^2}{M},
\end{eqnarray*}
   which yields the estimation in \eqref{3.06}.
\end{proof}

     By Lemma \ref{lem1} again, we can further divide the set $\{ 1, \ldots, k_2\}$ into the union of two subsets $I, J$ such that each $m$ can be written as
$\prod_{i \in I}p^{\alpha_i}_i\prod_{i \in J}p^{\alpha_i-1}_im'$.
There are $O_n(1)$ such partitions and without loss of generality,
we now consider the case
$m=\prod^{k_3}_{i=1}p^{\alpha_i}_i\prod^{k_2}_{i=k_3+1}p^{\alpha_i-1}_im'$,
for some integer $0 \leq k_3 \leq k_2$
   with $(m', \prod^{k_2}_{i=k_3+1}p_i)=1$ and here we take empty product to be $1$.
   We now bound $G_m(p^{\alpha_i}_i) \ll p^{\alpha_i}_i$ for $1 \leq i \leq k_3$ and
   $G_m(p^{\alpha_i}_i) \ll p^{\alpha_i-1}_i$ for $k_3+1 \leq i \leq k_2$ by Lemma
   \ref{lem1} and apply \eqref{3.05} with
$M=\prod^{k_3}_{i=1}p^{\alpha_i}_i\prod^{k_2}_{i=k_3+1}p^{\alpha_i-1}_i$
and $N=2\alpha^2Y/(\prod^{k_1}_{i=k_2+1}p^{\alpha_i-1}_i
   X)$ to bound the sum over those $m$'s in $R_1$ as
\begin{equation*}
 \ll  \frac {U
\alpha^2Y
\log^{2n+1}X}{\prod^{k_3}_{i=1}p^{\alpha_i}_i\prod^{k_1}_{i=k_3+1}p^{\alpha_i-1}_i
X}.
\end{equation*}

  Further note that $\alpha_i \geq 2$ for $1 \leq i \leq k_1$,
hence we deduce that
\begin{equation*}
  R_1 \ll \frac {U \alpha^2\sqrt{Y}\log^{2n+1} X}{X}.
\end{equation*}

  Similarly, by using the bound $G_m(p^{\alpha_i}_i) \ll p^{\alpha_i}_i$ for $1 \leq i \leq k_3$ and
   $G_m(p^{\alpha_i}_i) \ll p^{\alpha_i-1}_i$ for $k_3+1 \leq i \leq k_2$ and
  apply \eqref{3.05} with
$M=\prod^{k_3}_{i=1}p^{\alpha_i}_i\prod^{k_2}_{i=k_3+1}p^{\alpha_i-1}_i$
and
$N=2\alpha^2\prod^{k_2}_{i=1}p^{\alpha_i}_i\prod^{k_1}_{i=k_2+1}p_iV/X$,
we can bound the sum over those $m$'s in $R_2$ as
\begin{equation*}
  \ll \frac {U\alpha^2V(\prod^{k_1}_{i=k_3+1}p_i)\log^{2n+1}X}{X}.
\end{equation*}
   From this and that $\alpha_i \geq 3$ for $k_2+1 \leq i \leq k_1$, we deduce that
\begin{equation*}
   R_2 \ll \frac {U \alpha^2\sqrt{Y}\log^{3n+1} X}{X}.
\end{equation*}

   Lastly, we bound the sum over those $m$'s in $R_3$ by applying
    \eqref{3.06} with $M=\prod^{k_3}_{i=1}p^{\alpha_i}_i\prod^{k_2}_{i=k_3+1}p^{\alpha_i-1}_i$
and $N=2\alpha^2\prod^{k_2}_{i=1}p^{\alpha_i}_i
    \prod^{k_1}_{i=k_2+1}p_iV/X$ to be
\begin{equation*}
  \ll \frac {U^2\alpha^4(\prod^{k_2}_{i=1}p^{\alpha_i}_i
    \prod^{k_1}_{i=k_2+1}p_i)(\prod^{k_1}_{i=k_3+1}p_i)V^2
    \log^{2n+1}X}{X^2},
\end{equation*}
   which yields
\begin{equation*}
  R_3  \ll   \frac {U^2 \alpha^2 \sqrt{Y}\log^{3n+1} X}{X}.
\end{equation*}

   Combining the estimations above, we obtain
\begin{equation*}
  R \ll \frac {U^2Z \sqrt{Y}\log^{5n+2} X}{X}.
\end{equation*}
\subsection{The Second Main Term}
\label{sec 3.5}
  Using the notations in Section \ref{sec 3.3}, the sums corresponding to the contribution of $m=\Box$ in $S_{M}(X,Y;\prod^n_{i=1}{\hat g}_i, \Phi)$
 are linear
combinations of the sums of the following form:
\begin{eqnarray*}
 && \sum_{\substack {\prod^{k}_{i=1}p^{\alpha_i}_i \leq Y \\ (2, \prod^k_{i=1}p_i)=1 \\  p_i \neq p_j}}
   \prod^{k}_{i=1}\frac {(\log p_i)^{\alpha_i}}{p^{3\alpha_i/2}_i}
   h_i(\frac {\log p_i}{\log X}) \\
 && \cdot \sum_{\substack { \alpha \leq Z \\ (\alpha,
2\prod^k_{i=1}p_i)=1}}\frac {\mu(\alpha)}{\alpha^2}
    \sum^{\infty}_{m=1}(-1)^mG_m(p^{\alpha_i}_i) 
        \tilde{\Phi}\genfrac(){1pt}{}{m^2 X}{2\alpha^2\prod^{k}_{i=1}p^{\alpha_i}_i}.
\end{eqnarray*}

   Without loss of generality, we focus first on the case where
   there exists an integer $k_1$ with $0 \leq k_1 \leq k$
  such that for
   $k_1 < i \leq k$, $\alpha_i \equiv 1 \hspace{0.05in} (\text{mod} \hspace{0.05in} 2)$
    and $\alpha_i \equiv 0 \hspace{0.05in} (\text{mod} \hspace{0.05in} 2)$ for $i \leq
   k_1$.
   By Lemma \ref{lem1},
   each $m^2$ can be written as $m^2=\prod^{k_1}_{i=1}p^{\alpha_i}_i\prod^{k}_{i=k_1+1}p^{\alpha_i-1}_i(m')^2$
and we take empty product to be $1$. Thus the sum we are interested
 now becomes, by a slightly change of notation,
\begin{eqnarray}
\label{6.1}
 &&  \sum_{\substack {\prod^{k}_{i=1}p^{\alpha_i}_i \leq Y \\ (2, \prod^k_{i=1}p_i)=1 \\  p_i \neq p_j}}
   \prod^{k_1}_{i=1}\frac {(\log p_i)^{\alpha_i}}{p^{3\alpha_i/2}_i}\phi(p^{\alpha_i}_i)
   h_i(\frac {\log p_i}{\log X})\prod^{k}_{i=k_1+1}\frac {(\log p_i)^{\alpha_i}}{p^{(\alpha_i+1)/2}_i}
   h_i(\frac {\log p_i}{\log X})   \\
&&
 \cdot \sum_{\substack { \alpha \leq Z \\ (\alpha,
2\prod^k_{i=1}p_i)=1}}\frac {\mu(\alpha)}{\alpha^2}
    \sum^{\infty}_{\substack {m=1 \\ (m, \prod^{k}_{k_1+1}p_i)=1}}(-1)^m\tilde{\Phi}
    (\frac {m^2 X}{2\alpha^2\prod^{k}_{i=k_1+1}p_i}) \nonumber \\
  &=& \sum_{\substack {\prod^{k}_{i=1}p^{\alpha_i}_i \leq
Y \\ (2, \prod^k_{i=1}p_i)=1 \\  p_i \neq p_j}}
   \prod^{k_1}_{i=1}\frac {(\log p_i)^{\alpha_i}}{p^{3\alpha_i/2}_i}
 \phi(p^{\alpha_i}_i)
   h_i(\frac {\log p_i}{\log X})\prod^{k}_{i=k_1+1}
 \frac {(\log p_i)^{\alpha_i}}{p^{(\alpha_i+1)/2}_i}
   h_i(\frac {\log p_i}{\log X})  \nonumber \\
&&
  \cdot \sum_{\substack { \alpha \leq Z \\ (\alpha,
2\prod^k_{i=1}p_i)=1}}\frac {\mu(\alpha)}{\alpha^2} \sum_{d
|\prod^{k}_{k_1+1}p_i}\mu (d)
    \sum^{\infty}_{m=1}(-1)^m\tilde{\Phi}
    (\frac {d^2m^2 X}{2\alpha^2\prod^{k}_{i=k_1+1}p_i}), \nonumber
\end{eqnarray}
  where we have expressed the condition $(m,\prod^{k}_{k_1+1}p_i)=1$ as $\sum_{d|(m,\prod^{k}_{k_1+1}p_i)}\mu (d)$.
  To treat the sum over $m$ above, we need the following lemma:
\begin{lemma}
\label{lemma3.6}
  For $y>0$,
\begin{equation*}
   \sum^{\infty}_{m=1}(-1)^m\tilde{\Phi}(\frac {m^2}{y^2})=-\frac
   {\hat{\Phi}(0)}{2}+O(\frac {U}{y}).
\end{equation*}
\end{lemma}
\begin{proof}
   For $y>0$, we write
\begin{equation*}
   \sum^{\infty}_{m=1}(-1)^m\tilde{\Phi}(\frac {m^2}{y^2})=\frac 1{2}
   \sum^{\infty}_{m=-\infty}(-1)^m\tilde{\Phi}(\frac {m^2}{y^2})
   -\frac 1{2}\hat{\Phi}(0)=\frac {S(y)-\hat{\Phi}(0)}{2}.
\end{equation*}

   Define
\begin{equation*}
\label{38}
   F=\tilde{\Phi}(x^2),
\end{equation*}
   then by Poisson summation,
\begin{eqnarray}
   S(y) &=& \sum^{\infty}_{m=-\infty}(-1)^mF(\frac {m}{y})=2\sum^{\infty}_{m=-\infty}F(\frac {2m}{y})-
\sum^{\infty}_{m=-\infty}F(\frac {m}{y}) \nonumber \\
  &=& 2\sum^{\infty}_{\nu=-\infty}\hat{F}(\frac {y\nu}{2})\frac {y}{2}-\sum^{\infty}_{\nu=-\infty}\hat{F}(y\nu)y
= y\sum^{\infty}_{\substack{\nu=-\infty \\ (\nu, 2)=1}}\hat{F}(\frac {y\nu}{2}) \nonumber \\
&=&
y\sum^{\infty}_{m=-\infty}\int^{\infty}_{-\infty}\tilde{\Phi}(x^2)e\big(-xy(2m+1)/2
\big )dx
\nonumber \\
&=& 4y\sum^{\infty}_{m=0}\int^{\infty}_{0}\tilde{\Phi}(x^2)\cos(\pi y (2m+1)x) dx\nonumber \\
   &=& 4y\sum^{\infty}_{m=0}\int^{\infty}_{0}\tilde{\Phi}(x^2)d\frac {\sin(\pi y(2m+1)x)}{\pi y (2m+1)} \nonumber \\
   &=& - \frac 1{\pi}\sum^{\infty}_{m=0}\frac
   1{2m+1}\int^{\infty}_{0}2x\tilde{\Phi}'(x^2)\sin((2m+1)\pi
   yx)dx \nonumber \\
   &=& \frac 2{\pi^2 y}\sum^{\infty}_{m=0}\frac
   1{(2m+1)^2}\int^{\infty}_{0}x\tilde{\Phi}'(x^2)d\cos((2m+1)\pi
   y x) \nonumber \\
   &=& \frac {-2}{\pi^2
   y}\int^{\infty}_{0}\left(\tilde{\Phi}'(x^2)+2x^2\tilde{\Phi}''(x^2) \right )
   \left (\sum^{\infty}_{m=0}\frac
   {\cos((2m+1) \pi y x)}{(2m+1)^2} \right )dx \nonumber \\
 &\ll&  \frac {1}{y}\int^{\infty}_{0}\left (|\tilde{\Phi}'(x^2)|+|x^2\tilde{\Phi}''(x^2)| \right )
   dx
\label{39}
    \ll  \frac {1}{y}(\int^{1}_{0}1+\int^{\infty}_{1}\frac {U}{x^2})dx \ll \frac {U}{y}. \nonumber
\end{eqnarray}
   The last step above follows from \eqref{16'} and \eqref{18} with
   $i=1,2, j=2$ and this completes the proof.
\end{proof}

   We now divide the set $\{k_1+1, \ldots, k \}$ into the union of two subsets $I,J$ and set $d=\prod_{i \in I}p_i$.
Here we take empty product to be $1$. When $\prod_{j \in J}p_i \leq
X\prod_{j \in I}p_i$, we have
\begin{eqnarray*}
 &&   \sum^{\infty}_{m=1}|\tilde{\Phi}
    (\frac {m^2 X\prod_{i \in I}p_i}{2\alpha^2\prod_{i \in J}p_i})| \\
  &\ll& \sum_{m^2 X\prod_{i \in I}p_i <  2\alpha^2\prod_{i \in J}p_i}1+
  \sum_{m^2 X\prod_{i \in I}p_i \geq 2\alpha^2\prod_{i \in J}p_i}\frac {\alpha^2\prod_{i \in J}p_i}{m^2 X\prod_{i \in I}p_i}  \\
&\ll& \frac {\alpha \sqrt{\prod_{i \in J}p_i}}{\sqrt{X\prod_{i \in
I}p_i}}.
\end{eqnarray*}
    Hence these terms contribute to the right-hand side of \eqref{6.1} of order
\begin{equation*}
  \ll \log Z \log^{n-1}X.
\end{equation*}

   When $\prod_{j \in J}p_i > X\prod_{j \in I}p_i$, by our discussion above,
\begin{equation}
\label{6.3}
   \sum^{\infty}_{m=1}(-1)^m\tilde{\Phi}
    (\frac {m^2 X\prod_{i \in I}p_i}{2\alpha^2\prod_{i \in J}p_i})=-\frac {\hat{\Phi}(0)}{2} + O(\frac {U\sqrt{X\prod_{i \in I}p_i}}{\alpha \sqrt{\prod_{i \in J}p_i}}).
\end{equation}
   The error term above contribute to the right-hand side of \eqref{6.1} of order
\begin{equation*}
  \ll  U \log^{n-1}X.
\end{equation*}

   It's also easy to see that if $\alpha_i > 2$ for $i \leq k_1$ or $\alpha_i>1$ for $k_1 <  i \leq k$
then the first term on the right-hand side of \eqref{6.3} will
contribute to the right-hand side of \eqref{6.1} of order
\begin{equation*}
  \ll  \log^{n-1}X.
\end{equation*}

   Thus the main contribution to the right-hand side of \eqref{6.1}
   comes from the case $\alpha_i = 2, i \leq k_1$, $\alpha_i=1, k_1 <  i \leq
   k$, which corresponds to a partition of $\{1, \ldots, n \}$ into two subsets
   $S, S^c$ with $S^c$ nonempty, $|S|$ even and a way of pairing up the elements of $S$. Hence by partial summation,
we can write the sums in $S_{M}(X,Y;\prod^n_{i=1}{\hat g}_i, \Phi)$
corresponding to the contribution of $m=\Box$ as
\begin{eqnarray*}
    S_{\Box} &=&  - \frac {2X\log^n X \hat{\Phi}(0)}{\pi^2}\sum_{\substack {S \subsetneq \{1, \ldots, n\}
\\ |S| \text { even}}}
 \left ( \sum_{(A;B)}\prod^{|S|/2}_{i=1}\int^{\infty}_{0}u_i{\hat g}_{a_i}(u_i){\hat g}_{b_i}(u_i)du_i \right )  \\
 && \cdot \sum_{I \subsetneq S^c}(-1)^{|I|}\prod_{i \in I}\int^{\infty}_{0} {\hat g}_i(u_i)du_i
 \int_{\substack{(\mathbb{R}_{\geq 0})^{I^c} \\ \sum_{i \in I}u_i \leq \sum_{i \in I^c}u_i-1}}\prod_{i \in I^c}{\hat g}_i(u_i)du_i   \\
  &&  +O \left(X(U+\log Z) \log^{n-1} X + \frac {X\log^n X}{Z} \right ), \nonumber
\end{eqnarray*}
   where $\sum_{(A;B)}$ is over all ways of pairing up the elements of $S$ and we note here removing the condition
$p_i \neq p_j$ on the right-hand side of \eqref{6.1} only introduces
an error term of $O(X\log^{n-1}X)$.
\subsection{The Asymptotic Expression}
\label{sec 3.6}
    Putting together what we obtain and note that $\hat{\Phi}(0)=1+O(1/U)$, we
    get
\begin{eqnarray*}
 &&  S(X,Y;\prod^n_{i=1}{\hat g}_i, \Phi)  \\
&=& \frac {4X\log^n X }{\pi^2} \left(\frac {1+(-1)^n}{2} \right
)\sum_{(A;B)}\prod^{n/2}_{i=1}
\int^{\infty}_0 u_i{\hat g}_{a_i}(u_i){\hat g}_{b_i}(u_i)du_i  \\
 &&- \frac {2X\log^n X}{\pi^2}\sum_{\substack {S \subsetneq \{1, \ldots, n\} \\ |S| \text { even}}}
 \left ( \sum_{(C;D)}\prod^{|S|/2}_{i=1}\int^{\infty}_{0}u_i{\hat g}_{c_i}(u_i){\hat g}_{d_i}(u_i)du_i \right )  \\
 &&  \cdot  \sum_{I \subsetneq S^c}(-1)^{|I|}
 \prod_{i \in I}\int^{\infty}_{0} {\hat g}_i(u_i)du_i
 \int_{\substack{(\mathbb{R}_{\geq 0})^{I^c} \\ \sum_{i \in I}u_i \leq \sum_{i \in I^c}u_i-1} }
 \prod_{i \in I^c}{\hat g}_i(u_i)du_i
  \\
&& +O \left( X(U+\log Z) \log^{n-1}X+ \frac {X \log^{n} X}{U}+\frac {X \log^{2n+1} X}{Z} \right)  \\
&& +O(\frac {U^2Z \sqrt{Y}\log^{5n+2} X}{X}),
\end{eqnarray*}
   where $\sum_{(A;B)}$ is over all ways of pairing up the elements of $\{1,\ldots, n \}$ and
$\sum_{(C;D)}$ is over all ways of pairing up the elements of $S$.

 We now apply similar argument to the new choice $\Phi_1(t)=1-\Phi(t)$ for $1 \leq t \leq 2$
   and $\Phi_1(t)=0$ otherwise and
   by taking $Y=X^{2-\epsilon}, \epsilon = \frac {10(n+1)\log \log X}{\log X},U=\log \log X, Z=\log^{n+2} X$ to obtain
\begin{eqnarray*}
&&  S(X,Y;\prod^n_{i=1}{\hat g}_i)  \\
&=& \frac {4X\log^n X }{\pi^2} \left(\frac {1+(-1)^n}{2} \right
)\sum_{(A;B)}\prod^{n/2}_{i=1}
\int^{\infty}_0 u_i{\hat g}_{a_i}(u_i){\hat g}_{b_i}(u_i)du_i  \\
 &&- \frac {2X\log^n X}{\pi^2}\sum_{\substack {S \subsetneq \{1, \cdots, n \} \\ |S| \text { even}}}
 \left ( \sum_{(C;D)}\prod^{|S|/2}_{i=1}\int^{\infty}_{0}u_i{\hat g}_{c_i}(u_i){\hat g}_{d_i}(u_i)du_i \right )  \\
 && \cdot \sum_{I \subsetneq S^c}(-1)^{|I|}\prod_{i \in I}\int^{\infty}_{0} {\hat g}_i(u_i)du_i
 \int_{\substack{(\mathbb{R}_{\geq 0})^{I^c} \\ \sum_{i \in I}u_i \leq \sum_{i \in I^c}u_i-1}}
 \prod_{i \in I^c}{\hat g}_i(u_i)du_i
   \\
&& +O \left( X (\log \log X)^2 \log^{n-1}X+ \frac {X \log^{n}
X}{\log \log X} \right ).
\end{eqnarray*}

   From the above we deduce easily that
\begin{eqnarray}
\label{9.1}
   && \lim_{X \rightarrow \infty}\frac {\pi^2 }{4 X \log^n X}S(X,Y;\prod^n_{i=1}{\hat g}_i)  \\
   &=& \left(\frac {1+(-1)^n}{2} \right )\sum_{(A;B)}\prod^{n/2}_{i=1}
\int^{\infty}_0 u_i{\hat g}_{a_i}(u_i){\hat g}_{b_i}(u_i)du_i \nonumber  \\
 &&- \frac {1}{2}\sum_{\substack {S \subsetneq \{1, \ldots, n \} \\ |S| \text { even}}}
 \left ( \sum_{(C;D)}\prod^{|S|/2}_{i=1}\int^{\infty}_{0}u_i{\hat g}_{c_i}(u_i){\hat g}_{d_i}(u_i)du_i \right )
  \nonumber  \\
 &&  \cdot \sum_{I \subsetneq S^c}(-1)^{|I|}\prod_{i \in I}\int^{\infty}_{0} {\hat g}_i(u_i)du_i
 \int_{\substack{(\mathbb{R}_{\geq 0})^{I^c} \\
 \sum_{i \in I}u_i \leq \sum_{i \in I^c}u_i-1}}\prod_{i \in I^c}{\hat g}_i(u_i)du_i . \nonumber
\end{eqnarray}
\subsection{Conclusion}
\label{sec 3.7}
   We continue the discussion from Section \ref{sec 2.3} here. To consider the contribution of the product of the $D_l$'s,
    we need the following lemma \cite[Claim 1]{R} :
\begin{lemma}
\label{lem10.1}
   Suppose $\prod^n_{i=1}\hat{f}_i(u_i)$ is supported in $\sum^n_{i=1}|u_i| \leq \alpha$. Then
$\prod^k_{j=1}\hat{F}_{l_j}(u_j)$ is supported in $\sum^k_{j=1}|u_j|
\leq \alpha$.
\end{lemma}

  Now by the above lemma and \eqref{9.1},  the product of the $D_l$'s
  considered in the paragraph below Proposition \ref{prop2.5} contributes a factor of
\begin{eqnarray*}
  && \sum_{S_2 \subseteq S}\left ( (\frac {-1}{2})^{|S^c_2|} \prod_{l \in S^c_2}\int_{{\mathbb R}}\hat{F}_l(u)du \right )
   \cdot (-2)^{|S_2|}  \\
&&   \cdot \left ( \left (\frac {1+(-1)^{|S_2|}}{2} \right
)\sum_{(A;B)}\prod^{|S_2|/2}_{i=1}
\int^{\infty}_0 u_i\hat{F}_{a_i}(u_i)\hat{F}_{b_i}(u_i)du_i  \right. \\
 && - \frac {1}{2}\sum_{\substack {S_3 \subsetneq S_2 \\ |S_3| \text { even}}}
  \left ( \sum_{(C;D)}\prod^{|S_3|/2}_{i=1}\int^{\infty}_{0}u_i\hat{F}_{c_i}(u_i)\hat{F}_{d_i}(u_i)du_i \right )
  \\
 &&  \left. \cdot \sum_{I \subsetneq S^c_3}(-1)^{|I|}
 \prod_{i \in I}\int^{\infty}_{0} \hat{F}_i(u_i)du_i
 \int_{\substack{(\mathbb{R}_{\geq 0})^{I^c} \\
 \sum_{i \in I}u_i \leq \sum_{i \in I^c}u_i-1}}\prod_{i \in I^c}\hat{F}_i(u_i)du_i \right ).
\end{eqnarray*}
  Here $\sum_{(A;B)}$ is over all ways of pairing up the elements of $S_2$ and
$\sum_{(C;D)}$ is over all ways of pairing up the elements of $S_3$.
From this we find that \eqref{8.1} tends, as $X \rightarrow \infty$,
to the right-hand side of \eqref{10.06}.

  To conclude the proof of Theorem \ref{thm2.1}, it is left to show that dropping the term $O(\log \log X / \log X)$ won't affect our discussion
above, and this is given by the following lemma \cite[Lemma 2]{R} :
\begin{lemma}
\label{lem10.2}
  Let $a_l(d)=\sum_{\gamma_{8d}}F_l(L\gamma^{(j)}_{8d})$, with $F_l(x)=\prod_{i \in F_l}f_i(x)$, then
\begin{equation*}
   \lim_{X \rightarrow \infty}\frac {\pi^2}{4X}\sum_{d \in D(X)}\prod^{\nu(\underline{F})}_{l=1}a_l(d) 
  =  \lim_{X \rightarrow \infty}\frac {\pi^2}{4X}\sum_{d \in D(X)}\prod^{\nu(\underline{F})}_{l=1}
  \left( a_l(d) +O \Big (\log \log X / \log X \Big ) \right ).
\end{equation*}
\end{lemma}

\section*{Acknowledgement}
\setcounter{equation}{0}
   The results in this paper are part of the author's doctoral thesis at the University of Michigan (2005). We would like to thank both Professor Hugh Montgomery and Professor Kannan Soundararajan for their constant encouragements. We would also like to thank Professor Kannan Soundararajan for suggesting the problem to the author and many helpful conversations that lead to the results of the paper. Part of the writing of the paper was carried out while the author was visiting the American Institute of Mathematics in Fall 2005. The author thanks the American Institute of Mathematics for its generous support and hospitality. The author is also very grateful to Professor Brian Conrey for his encouragements during the author's visit at AIM and to Professor David Farmer for his NSF grant DMS-0244660 which made the author's visit to AIM possible.



\end{document}